\newcommand{\R}{\mathbb R}

\newcommand{\diver}{{{\rm{div}}}}

\newcommand{\grad}{{{\rm{grad}}}}

\documentclass[twoside,a4paper,12pt, draft]{amsart}

\usepackage{times}
\usepackage{dsfont}

\numberwithin{equation}{section}

\title[Riemannian submersions with discrete spectrum]{Riemannian submersions with discrete spectrum}
\author{G. Pacelli Bessa \and  J. Fabio Montenegro \and Paolo Piccione}
\date{January 5th, 2010}
\address{G. Pacelli Bessa \hfill\break\indent Departamento de Matem\'{a}tica, \hfill\break\indent
Universidade Federal do Ceara \hfill\break\indent Bloco 914 -- Campus do Pici\hfill\break\indent 60455-760 -- Fortaleza - CE -- Brazil \hfill\break\indent bessa@mat.ufc.br}
\address{J. Fabio Montenegro \hfill\break\indent Departamento de Matem\'{a}tica, \hfill\break\indent
Universidade Federal do Ceara \hfill\break\indent Bloco 914 -- Campus do Pici\hfill\break\indent 60455-760 -- Fortaleza - CE -- Brazil \hfill\break\indent fabio@mat.ufc.br}
\address{ Paolo Piccione\hfill\break\indent Departamento de Matem\'{a}tica-IME, \hfill\break\indent
Universidade de S\~{a}o Paulo \hfill\break\indent Rua do Mat\~{a}o 1010 \hfill\break\indent 05508-090 S\~{a}o Paulo-SP -- Brazil \hfill\break\indent piccione.p@gmail.com}
\thanks{The authors were partially supported by CNPq-CAPES (Brazil) and MEC project PCI2006-A7-0532 (Spain).}

\begin{document}

\theoremstyle{plain}\newtheorem*{theorem}{Theorem}
\theoremstyle{definition}\newtheorem*{defin*}{Definition}
\theoremstyle{plain}\newtheorem{teo}{Theorem}[section]
\theoremstyle{plain}\newtheorem{prop}[teo]{Proposition}
\theoremstyle{plain}\newtheorem{lem}[teo]{Lemma}
\theoremstyle{plain}\newtheorem*{lem-n}{Lemma}
\theoremstyle{plain}\newtheorem{cor}[teo]{Corollary}
\theoremstyle{definition}\newtheorem{defin}[teo]{Definition}
\theoremstyle{remark}\newtheorem{rem}[teo]{Remark}
\theoremstyle{plain} \newtheorem{assum}[teo]{Assumption}
\theoremstyle{remark}\newtheorem*{example-n}{Example}

\theoremstyle{plain}\newtheorem{teoabsolute}{Theorem}
\swapnumbers
\theoremstyle{definition}\newtheorem{example}{Example}

\maketitle
\begin{abstract}
We prove some estimates on the spectrum of the Laplacian of the total space of a Riemannian
submersion in terms of the spectrum of the Laplacian of the base and the geometry of the fibers.
When the fibers of the submersions are compact and  minimal, we prove that the total space is
discrete if and only if the base is discrete. When the fibers are not minimal, we prove a discreteness
criterion for the total space in terms of the relative growth of the mean curvature of the fibers
and the mean curvature of the geodesic spheres in the base.
We discuss in particular the case of warped products.
\end{abstract}

\section{Introduction}Let $M$ be a complete Riemannian manifold and  $\triangle\!=\!\diver \circ \grad $  be the Laplace-Beltrami operator
  acting on    the space of  smooth functions on $M$ with compact support.
 The operator $\triangle$ is  essentially self-adjoint, thus it has a unique  self-adjoint extension,
 to an unbounded operator, denoted by $\triangle$, whose domain is the set of  functions $ f\in L^{2}(M)$ so that $\triangle\! f\in L^{2}(M)$. Recall that
  the  spectrum of a self-adjoint operator $A$,  denoted by $\sigma(A)$,    is  formed by all $\lambda\in \mathbb{R}$ for which $A - \lambda I  $ is not injective or the inverse operator $(A- \lambda I)^{-1}$ is unbounded,  \cite{davies}.
  In this paper we are going to study the spectrum of  $-\triangle$, (the operator $\triangle$ is negative),  and  we refer to $\sigma(-\triangle)$ as the spectrum of $M$ and in this case only, we denote by $\sigma(M)$.   It is important (in our study) to distinguish the various types of  elements of the spectrum of $M$ in order to have a better understanding of the relations between $M$ and $\sigma (M)$. This way, it is said that the  set of all eigenvalues of $\sigma (M)$ is  the {\em point spectrum} $\sigma_{p}(M)$, while the {\em discrete spectrum} $\sigma_{d}(M)$ is the set of all isolated\footnote{Isolated in the sense that for some $\varepsilon>0$ one has that  $(\lambda-\varepsilon, \lambda + \varepsilon)\cap \sigma (\triangle^{M})=\lambda$.} eigenvalues  of finite multiplicity.  The {\em essential spectrum}   $\sigma_{\mathrm{ess}}(M)=\sigma (M)\setminus \sigma_{d}(M)$ is the complement of the discrete spectrum.

 There is a vast literature studying  the spectrum of complete Riemannian manifolds, among that, we point out  geometric restrictions  implying  that the spectrum   is purely  continuous ($\sigma_{p}(M)=\emptyset$), see \cite{donnelly1}, \cite{donnelly-garofalo}, \cite{escobar}, \cite{karp}, \cite{rellich}, \cite{tayoshi} or implying that the spectrum if discrete
 ($\sigma_{\mathrm{ess}}(M)\!=\!\emptyset$), see  \cite{baider}, \cite{bessa-jorge-montenegro}, \cite{donnelly-li}, \cite{harmer}, \cite{kleine1}, \cite{kleine2}.

In \cite{baider}, Baider  studied the essential spectrum of warped product manifolds
$W=X\times_{\gamma} Y=(X\times Y, dX^2+\gamma^{2}(x)dY^2)$,   where  $\gamma\colon X\to \mathbb{R}$
is a positive smooth function. The Laplace-Beltrami operator $\triangle^{W}$ restricted to
$ C_{0}^{\infty}(X)\otimes C_{0}^{\infty}(Y)$ has this form $\triangle^W\!\!= A_{0}\otimes 1_{Y}+
\gamma^{-2}\otimes (\triangle^{Y})$, where  $A_{0}$ is an elliptic operator, symmetric relative to
the density $\gamma^{n}dX$  with the same symbol as  $\triangle^X \!$.
 Baider   showed that if $\sigma_{\mathrm{ess}}(Y)=\emptyset$ then
 \begin{equation}\label{eqBaider}
 \sigma_{\mathrm{ess}}(W)=\emptyset\quad\Longleftrightarrow\quad \sigma_{\mathrm{ess}}(A_1)=\emptyset,
 \end{equation}
 where $A_{1}=A_{0}+ \lambda^{\ast}(Y)\gamma^{-2}$ and $\lambda^{\ast}(Y)=\inf \sigma (\triangle^{Y})$.   When $\gamma\equiv 1$ then $W=X\times Y$ and $A_{1}=\triangle^{X}+ \lambda^{\ast}(Y)$. In general, one can not  substitute $\sigma_{\mathrm{ess}}(A_1)=\emptyset$ by
 $\sigma_{\mathrm{ess}}(X)=\emptyset$. There are examples of warped manifolds $\mathbb{R}^{n}\times_{\gamma} \mathbb{S}^{1}\!,$ with discrete spectrum, therefore
 $\sigma_{\mathrm{ess}}(A_{1})=\emptyset$, see \cite{baider}, but the
 spectrum $\sigma_{\mathrm{ess}} (\mathbb{R}^{n})=[0, \infty)$. Riemannian manifolds
 whose Laplacian  has empty essential spectrum are sometimes called \emph{discrete} in the literature.

In this paper we consider Riemannian submersions $\pi \colon \! M\to N$ and we prove some spectral
estimates relating the (essential) spectrum of $M$ and $N$.
Riemannian submersions were introduced in the sixties by B. O'Neill and A. Gray (see \cite{Gra67,One66,One67})
as a tool to study the geometry of a Riemannian manifold with an additional structure in terms of certain
components, that is, the fibers and the base space.
When $M$ (and thus also $N$) is compact, estimates on the eigenvalues of the Laplacian of
$M$ have been studied in \cite{bordoni}, under the assumption that the mean curvature vector of the
fibers is \emph{basic}, i.e., $\pi$-related to some vector field on the basis.
We will consider here the non compact case, assuming initially
that the fibers are minimal.

An important class of examples are Riemannian homogeneous spaces $G/K$, where $G$ is a Lie group endowed with a bi-invariant
Riemannian metric and $K$ is a closed subgroup of $G$, see \cite{One66} for details. The projection $G\to G/K$ is a
Riemannian submersions with totally geodesic fibers, and with fibers diffeomorphic to $K$.

Another important class of examples of manifolds that can be described as the total space
of  Riemannian submersions with minimal fibers are  the homogeneous $3$-dimensional
Riemannian manifolds with isometry group of dimension four, see \cite{Sco}.
This class includes the special linear group $\mathrm{SL}(2,\R)$ endowed with a family of left-invariant
 metrics, which is the total space of  Riemannian submersions with base given by the hyperbolic spaces,
and fibers diffeomorphic to $\mathbb{ S}^1$.

Given a Riemannian submersion $\pi\!\colon \!\!M\!\to\! N$ with compact minimal fibers, we prove that
 \[\sigma_{\mathrm{ess}}(M)=\emptyset\quad \Longleftrightarrow\quad \sigma_{\mathrm{ess}}(N)=\emptyset,\] see Theorem~\ref{mainThm}.
This result coincides with Baider's result when $M=X\times Y$ is a product manifold,  $Y$ is compact, $N=X$ and $\pi\colon X\times Y\to X$ is the projection on the first factor.
\begin{teoabsolute}\label{mainThm}
Let $\pi\colon M\to N$ be a Riemannian submersion with compact
minimal fibers. Then \begin{itemize}\item[i.]$\sigma_{\mathrm{ess}}(N)\subset \sigma_{\mathrm{ess}}(M)$, $\sigma_{p}(N)\subset \sigma_{p}(M)$, thus  $\sigma(N)\subset \sigma (M)$.
\item[ii.]$\inf \sigma_{\mathrm{ess}}(N)=\inf \sigma_{\mathrm{ess}}(M)$. Therefore, $M$ is discrete if and only if  $N$
is discrete.
\end{itemize}
\end{teoabsolute}
A few remarks on this result are in order.
First, we observe that for the inequality $\inf\sigma_{\mathrm{ess}}(M)\leq
\inf\sigma_{\mathrm{ess}}(N)$, Lemma~\ref{lemma3.5}, we  need only the compactness of the fibers
with uniformly bounded volume, meaning that $0<c^{2}\leq {\rm vol}(\mathcal{ F}_{p})\leq C^{2} $ for all $p\in N$.
Second, the example of \cite{baider} shows that  the assumption of minimality of
the fibers is necessary in Theorem~\ref{mainThm}. In fact, one has examples of Riemannian submersions
having compact fibers with discrete base and non discrete total space, or with discrete total space
but not discrete base, see Example~\ref{examples}.

In the second part of the paper we study the essential spectrum of the total space when the
minimality assumption on the fiber is dropped. In this case, we prove that a sufficient condition
for the discreteness of the total space is that the growth of the mean curvature of the fibers at
infinity is controlled by the growth of the mean curvature of the geodesic spheres in the base manifold.
In order to state our result,  let us introduce the following terminology. The cut locus $cut(p)$ of
a point $p$ in a Riemannian $n$-manifold is said to be {\em thin}, if its $(n-1)$-Hausdorff measure zero,
$\mathcal{H}^{n-1}(cut (p))=0$.
\begin{teoabsolute}\label{mainThm2}
Let $\pi\colon M\to N$ be a Riemannian submersion with compact
fibers, and assume that $N$ has a point $x_0$ with thin cut locus. If
the function $h:M\to\R$ defined by
\[h(q)=(\triangle^N\rho_{p_0})_{\pi(q)}+g^N\big((\grad^N\rho_{p_0})_{\pi(q)},\mathrm d\pi_q(H_q)\big)\]
is \emph{proper} then
 $\sigma_{\mathrm{ess}}(M)=\emptyset$. Here $\rho_{p_0}$ is
the distance function in $N$ to $p_{0}$.
\end{teoabsolute}
The Theorem~\ref{mainThm2} can be interpreted geometrically in terms
of the mean curvature of the geodesic spheres in the base  and the mean curvature of the fibers.
Namely, the Laplacian of the distance function $\rho_{p_0}(p)$ is exactly
the value of the mean curvature of the geodesic sphere $\mathcal S_p=\rho_{p_0}^{-1}\big(\rho_{p_0}(p)\big)$
at the point $p$. Thus, assumption  says that the sum of the mean curvature of the geodesic balls in $N$
and the mean curvature of the fibers must diverge at infinity.

Theorem~\ref{mainThm} is proved in Section~\ref{sec:spectralestimate} and
Theorem~\ref{mainThm2} in Section~\ref{sec:mencurvfibersvscurvspheres}.
An alternative statement of Theorem~\ref{mainThm2} can be given in terms of radial curvature, see
Corollary~\ref{thm:radialcurvature}. There are two basic ingredients for the proof of our results.
\begin{itemize}
\item \noindent The Decomposition Principle, that relates the fundamental tone
of the complement of compact sets with the infimum of the essential spectrum, see Proposition~\ref{thm:charinfessspectrum};
\item Two estimates of the fundamental tones of open sets in terms of the divergence of vector fields, proved recently
 in \cite{bessa-montenegro} and \cite{bessa-montenegro2}, see Propositions~\ref{thm:est1} and \ref{thm:est2}.
\end{itemize}

\begin{section}{Riemannian submersions}
\subsection{Preliminaries}\label{sub:notterm}
Given  manifolds $M$ and $N$, a smooth surjective map $\pi\colon M\to N$   is a submersion
if the differential $\mathrm d\pi (q)$ has maximal rank for every $q\in M$.
If $\pi:M\to N$ is a submersion, then for all
$p\in N$ the inverse image $\mathcal F_p=\pi^{-1}(p)$ is a smooth embedded submanifold of $M$, that
will be called the \emph{fiber} at $p$. If $M$ and $N$ are Riemannian manifolds,
then a submersion $\pi:M\to N$ is called a \emph{Riemannian submersion} if for all $p\in N$ and
all $q\in\mathcal F_p$, the restriction of $\mathrm d\pi(q)$ to the orthogonal subspace
$T_q\mathcal F_p^\perp$ is an isometry onto $T_pM$.

Given $p\in N$ and $q\in\mathcal F_p$, a tangent vector $\xi\in T_qM$ is said to be \emph{vertical} if it
is tangent to $\mathcal F_p$, and it is \emph{horizontal} if it belongs to the orthogonal space
$(T_q\mathcal F_p)^\perp$. Let $\mathcal D=(T\mathcal F)^\perp\subset TM$ denote the smooth rank $k$ distribution
on $M$ consisting of horizontal vectors. The orthogonal distribution $\mathcal D^\perp$ is clearly
integrable, the fibers of the submersion being its maximal integral leaves.
Given $\xi\in TM$, its horizontal and vertical components are denoted respectively
by $\xi^{h}$ and $\xi^{v}$. The second fundamental form of the fibers is a symmetric tensor $\mathcal S^\mathcal F:\mathcal D^\perp\times\mathcal D^\perp\to\mathcal D$, defined by
\[\mathcal S^\mathcal F(v,w)=(\displaystyle\nabla^{M}_vW)^{h},\]
where  $W$ is a vertical extension of $w$ and $\nabla^M$ is the Levi--Civita connection of $M$.

For any given  vector field $X\in\mathfrak X(N)$,  there exists a unique horizontal
$\widetilde{X}\in\mathfrak X(M)$  which is $\pi$-related
to $X$, this is,  for any $p\in N$ and $q\in\mathcal F_p$,
then $\mathrm{ d}\pi_{q}(\widetilde{X}_{q})=X_{p}$, called {\em horizontal lifting} of $X$.   A horizontal vector field $\widetilde{X}\in\mathfrak X(M)$ is called \emph{basic} if it is $\pi$-related to
some vector field $X\in\mathfrak X(N)$.

\smallskip
If $\widetilde{X}$ and $\widetilde{Y}$ are basic vector fields, then these observations follows easily.
\begin{itemize}
\item[(a)] $g^M(\widetilde{X},\widetilde{Y})=g^N(X,Y)\circ\pi$.
\item[]
\item[(b)] $[\widetilde{X},\widetilde{Y}]^{h}$ is basic and it is $\pi$-related to $[X,Y]$.
\item[]
\item[(c)] $(\nabla^M_{\widetilde{X}}\widetilde{Y})^{h}$ is basic and it is $\pi$-related to $\nabla^N_{X}Y$,
\end{itemize}
where $\nabla^N$ is the Levi-Civita connection of $g^N$.

\smallskip
Let us now consider the geometry of the fibers. First, we observe that the fibers are totally geodesic
submanifolds of $M$ exactly when $\mathcal S^\mathcal F=0$.
The \emph{mean curvature} vector of the fiber is the horizontal vector field $H$ defined by
\begin{equation}H(q) = \sum_{i=1}^{k}\mathcal S^\mathcal F(q)(e_{i}, e_{i})=\sum_{i=1}^{k}(\displaystyle\nabla^{M}_{e_{i}}e_{i})^{h}\label{eq:defmeancurvature}
\end{equation}
where  $(e_i)_{i=1}^k$ is a local orthonormal frame for the fiber
through $q$. Observe that $H$ is not basic in general. For instance, when $n=1$, i.e., when the fibers
are hypersurfaces of $M$, then $H$ is basic if and only if all the fibers have constant
mean curvature. The fibers are \emph{minimal} submanifolds of $M$ when $H\equiv0$.
\subsection{Differential operators}

 Let  $\pi\colon \!\!M\to N$ be  a Riemannian submersion. Besides the  natural operations of lifting a
vector or vector fields in $N$ to horizontal vectors and basic vector fields one has that functions
on $N$ can be lifted to functions on $M$ that are constant along the fibers. Such operations
preserves the regularity of the lifted objects.  One can also (locally) lift  curves in the base $\gamma:[a,b]\to N$  to horizontal
curves $\widetilde{\gamma}:[a,c)\to M$ with the same regularity as $\gamma$
with arbitrary initial condition on the fiber $\mathcal F_{\gamma(a)}$.
We will need formulas relating the derivatives of  $\pi$-related objects in $M$ and $N$.
Let us start with divergence of vector fields.
\begin{lem}\label{thm:divergence}
Let $\widetilde{X}\!\in\!\mathfrak X(M)$ be a basic vector field, $\pi$-related to $X\!\in\!\mathfrak X(N)$.
The following relation holds between the divergence of $\widetilde{X}$ and $X$ at $p\!\in\!N$ and $q\!\in\!\mathcal F_p$.
\begin{equation}\label{eq:reldivergence}\begin{array}{lll}
\mathrm{div}^M(\widetilde{X})_q&=&\mathrm{div}^N(X)_p+g^M(\widetilde{X}_q,H_q) \\
&& \\ &=&\mathrm{div}^N(X)_p+g^N\big(\mathrm d\pi_q(\widetilde{X}_q),\mathrm d\pi_q(H_q)\big).\end{array}
\end{equation}
In particular, if the fibers are minimal, then $\mathrm{div}^M(\widetilde{X})=\mathrm{div}^N(X)$.
\end{lem}
\begin{proof}
Formula \eqref{eq:reldivergence} is obtained by a direct computation of the left-hand side,
using a local orthonormal frame $e_1,\ldots,e_k,e_{k+1},\ldots,e_{k+n}$ of $TM$, where
$e_1,\ldots,e_k$ are basic fields. The  equality follows  using equalities (a) and (c) in
Subsection~\ref{sub:notterm}, and formula \eqref{eq:defmeancurvature} for the mean curvature.
\end{proof}
Given a smooth function $f:N\to\R$, denote by $\tilde f=f\circ\pi:M\to\R$ its \emph{lifting} to $M$.
It is easy to see that the gradient $\grad^M\tilde f$ of $\tilde f$ is the horizontal lifting of
the gradient $\grad^Nf$. If we denote with a tilde $\widetilde X$ the horizontal lifting of a vector
field $X\in\mathfrak X(N)$, then the previous statement can be written as
\begin{equation}\label{eq:gradlifting}
\grad^M\tilde f=\widetilde{\grad^Nf}.
\end{equation}
Now, given a function $f\!:\!M\to\mathbb{R}$, one can define a function $f_{\mathrm{av}}\!:\!N\to\mathbb{R}$ by averaging
$f$ on each fiber
\[f_{\mathrm{av}}(p)=\int_{\mathcal F_p}\!\!f\,\mathrm d\mathcal F_p,\]
where $\mathrm d\mathcal F_p$ is the volume element of the fiber $\mathcal F_p$ relative to the induced
metric. We are assuming that this integral is finite. As to the gradient of the averaged function $f_{\mathrm{av}}$, we have the following lemma.
\begin{lem}
Let $p\in N$ and $v\in T_pN$ and
denote by $V$ the smooth normal vector field along $\mathcal F_p$ defined by the property
$\mathrm d\pi_q(V_q)=v$ for all $q\in\mathcal F_p$.
Then, for any smooth function $f:M\to\R$
\begin{equation}\label{eq:gradaverage}
g^N\big(\grad^Nf_{\mathrm{av}}(p),v\big)\!=\!\!\int_{\mathcal F_p}\!\left[g^M\big(\grad^Mf,V\big)+f\!\cdot\! g^M(H,V)\right]\mathrm  d\mathcal F_q.
\end{equation}
\end{lem}
\begin{proof}
A standard calculation as in the first variation formula for the volume functional of the fibers.
Notice that when $f\equiv1$, then $f_{\mathrm{av}}$ is the volume function of the fibers,
 and
\eqref{eq:gradaverage} reproduces the first variation formula for the volume.
\end{proof}
Observe that, in \eqref{eq:gradaverage}, the gradient $\grad^Mf$ need not be basic or even horizontal\footnote{%
In fact, a gradient is basic if and only if it is horizontal.}.
An averaging procedure is available also to produce vector fields $X_{\mathrm{av}}$ on the basis out of vector fields $X$
defined in the total space. If $X\in\mathfrak X(M)$, let $X_{\mathrm{av}}\in\mathfrak X(N)$ be defined by
\[(X_{\mathrm{av}})_p=\int_{\mathcal F_p}\mathrm d\pi_{q}\big(X_q\big)\,\mathrm d\mathcal F_p(q).\]
Observe that the integrand above is a function on $\mathcal F_p$ taking values in the fixed vector space
$T_pN$. If $X\in\mathfrak X(M)$ is a basic vector field, $\pi$-related to the vector field $X_{*}\in\mathfrak X(N)$,
then $(X_{\mathrm{av}})_p=\mathrm{vol}(\mathcal F_p)\cdot(X_{*})_p$, where $\mathrm{vol}$ denotes the volume.
Using the notion of averaged field, equality \eqref{eq:gradaverage} can be rewritten as
\[\grad^N(f_{\mathrm{av}})=\big(\grad^Mf+f\cdot H\big)_{\mathrm{av}}.\]
\begin{rem}\label{thm:remminimalfibers}
From the above formula it follows easily that the averaged mean
curvature vector field $H_{\mathrm{av}}$ vanishes at the point $p\in N$ if and only if $p$ is a critical
point of the function $z\mapsto\mathrm{vol}(\mathcal F_z)$ in $N$. This happens, in particular,
when the leaf $\mathcal F_p$ is minimal. When all the fibers are minimal, or more generally
when the averaged mean curvature vector field $H_{\mathrm{av}}$ vanishes identically,
then the volume of the fibers is constant.
\end{rem}
\begin{cor}Let $\pi\colon M\to N$ be a Riemannian submersion with compact minimal fibers $\mathcal{F}$. Let $h\in L^{2}(N)$. If  $f \in C^{\infty}_{0}(M)$ such that $f_{\mathrm{av}}=0$ for all $q\in N$ then \begin{equation}\label{eq mediazero} \int_{M} \widetilde{h}\,\triangle^{M}f \, dM=0.
\end{equation}\label{thmmediazero}
\end{cor}
\begin{proof} Suppose first that $h$ is smooth. By the Divergence Theorem, Fubini's Theorem for Riemannian submersions and \ref{eq:gradaverage} we have  \begin{eqnarray}\int_{M} \widetilde{h}\,\triangle^{M}f \, dM& =& -\int_{M}g^{M}(\grad ^{M}\widetilde{h}, \grad ^{M} \!f)dM\nonumber \\
&=& -\int_{N}\int_{\mathcal{F}_{q}}g^{M}(\grad ^{M}\widetilde{h}, \grad ^{M} \!f)d\mathcal{F}_{q}dN\nonumber\\ &=& -\int_{N}g^{N}(\grad ^{M}\widetilde{h}, \grad ^{N} \!f_{\mathrm{av}})dN\nonumber\\ &=& 0 \nonumber
\end{eqnarray} If  $h\in L^{2}(N)$  there exists a sequence of smooth functions $h_{k}\in C^{\infty}(N)$ converging to $h$ with respect to the $L^{2}$-norm.
On the other hand
\begin{eqnarray}\left \vert \int_{M}\widetilde{h}\triangle^{M}f \,dM \right\vert & = &
\left \vert \int_{M}(\widetilde{h}_{k}-\widetilde{h})\triangle^{M}f \,dM \right\vert\nonumber \\ & \leq & \int_{M}\left\vert \widetilde{h}_{k}-\widetilde{h}\right\vert \, \vert \triangle^{M}f \vert \, dM \nonumber \\
&\leq & \!\left(\int_{M}\big\vert \widetilde{h}_{k}-\widetilde{h}\big\vert^{2}\!dM\right)^{1/2}\!\! \cdot\left(\int_{M}\!\vert \triangle^{M}f \vert^{2}  dM\right)^{1/2}\nonumber \\
&=&\!\Vert \triangle^{M}f\Vert_{L^{2}(M)}\cdot\left( \int_{N}\int_{\mathcal{F}_{q}}\big\vert h_{k}-h\big\vert^{2}d\mathcal{F}_{q}dM\right)^{1/2}\nonumber \\
&=& \!\mathrm{vol}(\mathcal{F}_{q})^{1/2}\cdot\Vert \triangle^{M}f\Vert_{L^{2}(M)}\cdot\Vert h_{k}-h\Vert_{L^{2}(N)}\nonumber
\end{eqnarray} Since $h_{k}\to h$ in $L^{2}(N)$ then \ref{eq mediazero} holds. Observe that we used that  the volume of the minimal fibers is constant, see Remark \ref{thm:remminimalfibers}.
\end{proof}
\end{section}

\begin{section}{Spectral estimates in Riemannian submersions}
\label{sec:spectralestimate}
\subsection{Generalities on  the Laplace-Beltrami operator}

Let $\Omega \subset M$ be an open set in complete Riemannian manifold. The fundamental tone of $\Omega$ is defined by
$$\lambda^{\ast}(\Omega)=\inf\frac{\int_{\Omega}\vert\grad f\vert^2}{\int_{\Omega}f^2},$$
where the infimum is taken over all smooth non zero functions $f\in C^{\infty}_{0}(\Omega)$.

The fundamental tone has the following monotonicity property:
if $A\!\subset \!B$ are open subsets then $\lambda^*(A)\ge\lambda^*(B)$. If $\Omega$ is an open subset of $M$ with compact closure, or if $\lambda^*(\Omega)\not\in\sigma_{\mathrm{ess}}(M)$,
then $\lambda^*(\Omega)$ coincides with the first eigenvalue $\lambda_1(\Omega)$
of the Dirichlet problem
$$ \left\{\begin{array}{rcl}\triangle^Mu+\lambda u\!\! \!\! & =&\!\! \! \! 0\,\,\text{in} \,\,\Omega,\\ u\!\!\!\!&=&\!\! \!  \!0\,\,\text{on} \,\,\partial\Omega.\end{array}\right.$$ When $\Omega =M$  then $\lambda^{\ast}(M)=\inf \sigma (M)$.

The Laplace-Beltrami operator in $N$ of a smooth function $f:N\to\R$ and the Laplace-Beltrami operator in $M$ of its extension
$\tilde f=f\circ\pi$ are related by the following formula.
\begin{lem}\label{thm:laplacianlift}
Let $f:N\to\R$ be a smooth function and set $\tilde f=f\circ\pi$. Then, for all $p\in N$ and all $q\in\mathcal F_{p}$:
\begin{equation}\label{eq:compLaplacianfbarf}\begin{array}{lll}
(\triangle^M\tilde f)_q&=&(\triangle^Nf)_p+g^M\big((\grad^M\tilde f)_q,H_q\big)\\ && \\ & = &
(\triangle^Nf)_p+g^N\big((\grad^Nf)_p,\mathrm d\pi_q(H_q)\big).\end{array}
\end{equation}
\end{lem}
\noindent The proof
 follows easily from \eqref{eq:reldivergence} applied to the vector fields $X=\grad^M\tilde f$ and $X_*=\grad^Nf$,
using \eqref{eq:gradlifting}.

\subsection{Decomposition Principle} Let $K\subset M$ be a compact  set of the same dimension of $M$.  The Laplace-Beltrami operator $\triangle$ of $M$ acting on the space $C^{\infty}_{0}(M\setminus K)$ of smooth compactly supported functions of $M\setminus K$  has a self-adjoint extension, denoted by $\triangle'$.
The {\em Decomposition Principle} \cite{donnelly-li} says that   $\sigma_{\mathrm{ess}}(M)=\sigma_{\mathrm{ess}}(M\setminus K).$
On the other hand,
\[0\le\lambda^*(M\setminus K)=\inf\sigma(M\setminus K)\le\inf\sigma_{\mathrm{ess}}(M\setminus K)=\sigma_{\mathrm{ess}}(M),\] thus $\mu = \sup\{\lambda^{\ast}(M\setminus K), K \subset M\,\,{\rm compact}\}\leq \sigma_{\mathrm{ess}}(M)$. We will show that $\inf \sigma_{\mathrm ess}(M)\leq \mu$. To that we will suppose that $\mu < \infty$, otherwise there is nothing to prove.
Let $K_{1}\subset K_{2}\subset \cdots $ be a sequence of compact sets with $M=\cup_{i=1}^{\infty}K_{i}$. We have that $$ \lambda^{\ast}(M)\leq \lambda^{\ast}(M\setminus K_{1})\leq \lambda^{\ast}(M\setminus K_{2})\leq \cdots  \to  \mu .$$
Given $\varepsilon >0$, there exists $f_{1}\in C_{0}^{\infty}(M\setminus K_{1})$ with $\Vert f_{1}\Vert_{L^{2}}=1$ and $\int_{M}\vert \grad f_{1}\vert^{2} \leq \lambda^{\ast}(M\setminus K_{1}) + \varepsilon <  \mu + \varepsilon $. This is
 $\langle (-\triangle -\mu -\varepsilon) f_{1}, f_{1}\rangle_{L^{2}}<0. $  We can suppose that $\mathrm{supp}\, (f_{1}) \subset (K_{2}\setminus K_{1})$. There exists $f_{2}\in C_{0}^{\infty}(M\setminus K_{2})$ with $\Vert f_{2}\Vert_{L^2}=1$ and $\int_{M}\vert \grad f_{2}\vert^{2} \leq \lambda^{\ast}(M\setminus K_{2}) + \varepsilon <  \mu + \varepsilon $. This is equivalent to  $\langle (-\triangle -\mu -\varepsilon) f_{2}, f_{2}\rangle_{L^{2}}<0. $
 Moreover, $\int_{M}f_{1}f_{2}=0$ since $\mathrm{supp}\, (f_{1})\cap \mathrm{supp}\, (f_{2})=\emptyset$.
 This way,   we obtain an orthonormal sequence $\{f_{k}\}\subset C_{0}^{\infty}(M)$ such that
 $\langle (-\triangle -\mu -\varepsilon) f_{k}, f_{k}\rangle_{L^{2}}<0$.
 By
\ref{lemmaPaolo} we have that $(-\infty, \mu]\cap \sigma_{\mathrm{ess}}(M)\neq \emptyset$ and $\inf \sigma_{\mathrm{ess}}(M)\leq \mu$. This proves the following proposition.
\begin{prop}\label{thm:charinfessspectrum}The infimum of the essential spectrum is characterized by
 \begin{equation}\label{eq:charinfessspectrum}
\inf\sigma_{\mathrm{ess}}(M)=\sup\big\{\lambda^*(M\setminus K):\ \text{$K$ compact subset of $M$}\big\}.
\end{equation}In particular, $\sigma_{\mathrm{ess}}(M)$ is empty if and only if given any compact
exhaustion $K_1\subset K_2\subset\cdots\subset K_n\subset\ldots$ of $M$, the
limit $\lim\limits_{n\to\infty}\lambda^*(M\setminus K_n)$ is infinite.
\end{prop}

 Let $H$ be a Hilbert space and $A\colon\mathcal D\subset\! H\to H$ be a densely defined self-adjoint  operator.
Given $\lambda\in\mathds R$, we write $A\ge\lambda$ if $\langle Ax,x\rangle\ge\lambda\Vert x\Vert^2$ for all $x\in\mathcal D$. By the Spectral Theorem for (unbounded) self-adjoint operators, we have that $A\ge\lambda$ iff $\sigma(A)\subset[\lambda,+\infty)$.
Let us write $A>-\infty$ if there exists $\lambda_*\in\mathds R$ such that $A\ge\lambda_*$.
\begin{lem}\label{lemmaPaolo}
Let $A\colon\mathcal D\subset \!H\to \!H$ be a self-adjoint operator with $A>-\infty$, and let $\lambda\in\mathds R$
be fixed. Assume that for all $\varepsilon>0$ there exists an infinite dimensional subspace
$G_\varepsilon\subset\mathcal D$ such that $\langle Ax,x\rangle<(\lambda+\varepsilon)\Vert x\Vert^2$
for all $x\in G_\varepsilon$. Then, \[\sigma_{\mathrm{ess}}(A)\cap(-\infty,\lambda]\ne\emptyset.\]
\end{lem}This lemma is well known, see  \cite{donnelly}  but for sake of completeness we present here its proof.
\begin{proof}
 First we  will show that  $\sigma(A)\cap(-\infty,\lambda]=\sigma(A)\cap[\lambda_{*}, \lambda]\ne\emptyset$.
Take $\varepsilon_k=1/k$, $k\ge1$. By our hypothesis there exists $x_k\ne0$
such that $\langle Ax_k,x_k\rangle<(\lambda+1/k)\Vert x_k\Vert^2$, and thus
$\sigma(A)\cap[\lambda_{*},\lambda+1/k]\ne\emptyset$ for all $k\ge1$.
Since $\sigma(A)$ is closed, it follows $\sigma(A)\cap(-\infty,\lambda]\ne\emptyset$.
We may suppose  that $\sigma(A)\cap(-\infty,\lambda]
\not\subset\sigma_{\mathrm{ess}}(A)$, otherwise there is nothing to prove. Thus \[\big(\sigma(A)\setminus\sigma_{\mathrm{ess}}(A)\big)
\cap(-\infty,\lambda]=\{\lambda_1,\ldots,\lambda_n\}\] is a finite set of eigenvalues
of $A$ of finite multiplicity.  Denote by $H_i\subset\mathcal D$
the $\lambda_i$-eigenspace of $A$, $i=1,\ldots,n$, and set $X=\bigoplus_iH_i\subset\mathcal D$.
This is clearly an invariant subspace of $A$. Since $X$ has finite dimension, then
$\mathcal D=X\oplus X_{1}$ where $X_1=X^\perp\cap\mathcal D$ is also invariant by $A$.
Denote by $A_1$ the restriction of $A$ to the Hilbert space $X_{1}$ which is still self-adjoint. Clearly, $\sigma(A_1)=\sigma(A)\setminus\{\lambda_1,\ldots,\lambda_n\}$
and $\sigma_{\mathrm{ess}}(A_1)=\sigma_{\mathrm{ess}}(A)$.  In particular, we have
$\sigma(A_1)\cap(-\infty,\lambda]\subset\sigma_{\mathrm{ess}}(A_1)$.
Using the infinite dimensionality of the space $G_\varepsilon$, it is now easy to see
that the assumptions of our lemma hold for the operator $A_1$, and the first part of the proof
applies to obtain $\sigma_{\mathrm{ess}}(A)\cap(-\infty,\lambda]=\sigma_{\mathrm{ess}}(A_1)
\cap(-\infty,\lambda]\ne\emptyset$.
\end{proof}
Let us  recall from \cite{bessa-montenegro} and \cite{bessa-montenegro2} the following estimates for the fundamental tone of open sets of Riemannian manifolds.
\begin{prop}\label{thm:est1}
Let $\Omega\subset M$ be an open set of a Riemannian manifold. Then
\begin{equation}\label{eq:secondestimate}
\lambda^*(\Omega)\ge\frac{1}{4}\sup_X\left[\frac{\inf\limits_{\Omega}\mathrm{div}(X)}{\sup\limits_{\Omega}\Vert X\Vert}\right]^2\!\!,
\end{equation}
where the supremum in taken over all smooth vector fields $X$ in $\Omega$ satisfying
\[\inf\limits_{\Omega}\mathrm{div}(X)>0,\quad\text{and}\quad\sup\limits_{\Omega}\Vert X\Vert<+\infty.\]
\end{prop}
\begin{prop}\label{thm:est2}
Let $\Omega\subset M$ be an open set of a Riemannian manifold.
Given any smooth vector field $X\in\mathfrak X(\Omega)$ then
\begin{equation}\label{eq:estimatefundtone}
\lambda^*(\Omega)\ge\inf_{\Omega}\left[\mathrm{div}(X)-\big\vert X\big\vert^2\right].
\end{equation}
Equality in \eqref{eq:estimatefundtone} holds when $\lambda^*(\Omega)\in\sigma(\Omega)\setminus\sigma_{\mathrm{ess}}(\Omega)$,
by considering the field $X=-\nabla(\log f)$, where $f$ is the positive eigenfunction associated to $\lambda^*(\Omega)$.
\end{prop}
\begin{rem}Propositions (\ref{thm:est1}) and (\ref{thm:est2}) hold  for
  vector fields  $X$ of class $ C^{1}(\Omega\setminus F)\cap
L^{\infty}(\Omega)$ such that $\diver (X)\in L^{1}(\Omega)$, where $F\subset M$ is a closed subset with
$(n-1)$-Hausdorff measure ${\mathcal H}^{n-1}(F\cap \Omega)=0$, see \cite[Lemma 3.1]{bessa-montenegro2}.

\label{thincutlocus}
\end{rem}
\subsection{Fundamental tones  estimates of Riemannian submersions}
Let $M$ and $N$ be connected Riemannian manifolds and $\pi\!\colon \!M\!\to\! N$ be a Riemannian submersion.
Denote by $\triangle^M$ and $\triangle^N$ the Laplacian operator
on functions of $(M,g^M)$ and of $(N,g^N)$ respectively. We want to compare the fundamental tones of open subsets $\Omega \subset N$ with the fundamental tones of its lifting $\widetilde{\Omega}=\pi^{-1}(\Omega)$.
\begin{lem}\label{lemma3.5}
Assume that the fibers of $\pi:M\to N$ are compact.
Let $\Omega$ be an open subset of $N$, and denote by $\widetilde\Omega$ the open subset of $M$
given by the inverse image $\pi^{-1}(\Omega)$.
Then
\begin{equation}\label{eq:estimate}
\left[{\inf\limits_{p\in\Omega}\mathrm{vol}(\mathcal F_p)}\right]\cdot\lambda^*(\widetilde\Omega)\le
\left[{\sup\limits_{p\in\Omega}\mathrm{vol}(\mathcal F_p)}\right]\cdot\lambda^*(\Omega).
\end{equation}

In particular, if the fibers are minimal, then
\begin{equation}\label{eq:betterestimate}\lambda^*(\widetilde\Omega)\le\lambda^*(\Omega).\end{equation}
Moreover, if $\inf\limits_{p\in\Omega}\mathrm{vol}(\mathcal F_p)>0$ and ${\sup\limits_{p\in\Omega}\mathrm{vol}(\mathcal F_p)}<\infty$ then
\begin{equation}\label{eq:betterestimate2}\inf\limits_{p\in\Omega}\mathrm{vol}(\mathcal F_p)\cdot\inf\sigma_{\mathrm{ess}}(M)\le {\sup\limits_{p\in\Omega}\mathrm{vol}(\mathcal F_p)}\cdot\inf\sigma_{\mathrm{ess}}(N).\end{equation}
\end{lem}
\begin{proof}
Let $\varepsilon>0$  and choose $f_\varepsilon\in C^\infty_0(\Omega)$ such that
\begin{equation}\label{eq:basicestfepsilon}
{\int_\Omega\big\vert\grad^Nf_\varepsilon\big\vert^2}<\big(\lambda^*(\Omega)+\varepsilon\big)
{\int_\Omega f_\varepsilon^2}.
\end{equation}
Let us consider the function $\tilde f_\varepsilon=f_\varepsilon\circ\pi$. By the assumption that
the fibers of $\pi$ are compact, $\tilde f_\varepsilon$ has compact support in $M$.
Using Fubini's Theorem for submersions we have
\begin{multline*}
\int_{\widetilde\Omega}\big\vert\tilde f_\varepsilon\big\vert^2\;\mathrm dM=\int_\Omega\left(\int_{\mathcal F_p}
\big\vert\tilde f_\varepsilon\big\vert^2\;\mathrm d\mathcal F_p\right)\mathrm dN
=\int_\Omega\mathrm{vol}(\mathcal F_p)\cdot\vert f_\varepsilon\big\vert^2\mathrm dN.
\end{multline*}
Thus
\begin{equation}\label{eq:L2normftilde}
\int_{\widetilde\Omega}\big\vert\tilde f_\varepsilon\big\vert^2\mathrm dM\ge
{\inf\limits_{p\in\Omega}\mathrm{vol}(\mathcal F_p)}
\cdot\int_{\Omega}\big\vert f_\varepsilon\big\vert^2\mathrm dN.
\end{equation}
Similarly, using \eqref{eq:gradlifting}, we have
\begin{eqnarray}
\int_{\widetilde\Omega}\big\vert\grad^M\tilde f_\varepsilon\big\vert^2 & = &
\int_{\widetilde\Omega}\big\vert\widetilde{\grad^Nf_\varepsilon}\big\vert^2 \nonumber \\ & = &
\int_\Omega\Big(\int_{\mathcal F_p}
\big\vert\widetilde{\grad^Nf_\varepsilon}\big\vert^2\;\mathrm d\mathcal F_p\Big)\mathrm dN\\ & = &
\int_\Omega\mathrm{vol}(\mathcal F_p)\cdot\big\vert{\grad^Nf_\varepsilon}\big\vert^2,\nonumber
\end{eqnarray}
thus
\begin{equation}\label{eq:L2normgradftilde}
\int_{\widetilde\Omega}\big\vert\grad^M\tilde f_\varepsilon\big\vert^2\le
{\sup\limits_{p\in\Omega}\mathrm{vol}(\mathcal F_p)}
\cdot\int_\Omega\big\vert{\grad^Nf_\varepsilon}\big\vert^2.
\end{equation}
Using \eqref{eq:basicestfepsilon}, \eqref{eq:L2normftilde} and \eqref{eq:L2normgradftilde}, we then
obtain
\begin{eqnarray}
{\inf\limits_{p\in\Omega}\mathrm{vol}(\mathcal F_p)}\cdot\lambda^*(\widetilde\Omega)
& \le & {\inf\limits_{p\in\Omega}\mathrm{vol}(\mathcal F_p)}\cdot
\frac{\int_{\widetilde\Omega}\big\vert\grad^M\tilde f_\varepsilon\big\vert^2}{\int_{\widetilde\Omega}\big\vert
\tilde f_\varepsilon\big\vert^2}\nonumber \\
& \le & {\sup\limits_{p\in\Omega}\mathrm{vol}(\mathcal F_p)}\cdot
\frac{\int_\Omega\left\vert\grad^Nf_\varepsilon\right\vert^2}{\int_\Omega\vert f_\varepsilon\vert^2} \\ && \nonumber \\ &< &
{\sup\limits_{p\in\Omega}\mathrm{vol}(\mathcal F_p)}\cdot\left[\lambda^*(\Omega)+\varepsilon\right].\nonumber
\end{eqnarray}
This proves (\ref{eq:estimate}). If all the fibers are minimal
(or more generally if the averaged mean curvature vector field $H_{\mathrm{av}}$ vanishes
identically on $N$, see Remark~\ref{thm:remminimalfibers}), then the volume of the fibers
is constant, and inequality \eqref{eq:betterestimate} follows from \eqref{eq:estimate}. To prove the inequality \eqref{eq:betterestimate2} we
 pick a compact subset $K\subset M$ and set $K_0=\pi(K)$ and let $\widetilde K=\pi^{-1}(K_0)$. The set $\widetilde K$
is compact by the assumption that the fibers of $\pi$ are compact.
Let $\Omega=N\setminus K_0$ and  $\widetilde\Omega=\pi^{-1}(\Omega)=M\setminus\widetilde K$.
Clearly, $\widetilde\Omega\subset M\setminus K$ and thus $\lambda^*(\widetilde\Omega)\ge\lambda^*(M\setminus K)$.
Hence, using \eqref{eq:estimate} we get
\[\lambda^*(M\setminus K)\le\lambda^*(\widetilde\Omega)\le\frac{\sup\limits_{p\in\Omega}\mathrm{vol}(\mathcal F_p)}{\inf\limits_{p\in\Omega}\mathrm{vol}(\mathcal F_p)}\,\lambda^*(\Omega)\le
\frac{\sup\limits_{p\in\Omega}\mathrm{vol}(\mathcal F_p)}{\inf\limits_{p\in\Omega}\mathrm{vol}(\mathcal F_p)}\,\inf\sigma_{\mathrm{ess}}(N).\]
Taking the supremum over all compact subset $K\subset M$ in the left-hand side, we obtain the desired inequality.
\end{proof}

Now   consider the case that the fibers of the submersion $\pi \colon M \to N$ are compact and minimal.
\begin{lem}\label{thm:equalfundtone}
Let  $\pi\!:\!M\!\to\! N$ be a Riemannian submersion with  compact and minimal fibers ${\mathcal F}$.
Then for every open subset $\Omega\subset N$, denoting by $\widetilde\Omega$ the inverse image
$\pi^{-1}(\Omega)$, one has that
\begin{equation}\label{eq:equalfundtone}
\lambda^*(\widetilde\Omega)=\lambda^*(\Omega).
\end{equation}
\end{lem}

\begin{proof}
In view of \eqref{eq:betterestimate}, it suffices to show the inequality $\lambda^*(\widetilde\Omega)\ge
\lambda^*(\Omega)$. To this aim, we will use the estimate in \eqref{eq:estimatefundtone}.
We observe initially that it suffices to prove the inequality when $\Omega$ is bounded.
Namely, the general case follows from $\lambda^{\ast}(\Omega)=\lim_{n\to \infty} \lambda^{\ast}(\Omega_{n})$, by considering an exhaustion
of $\Omega$ by a sequence of bounded open subsets $\Omega_{n}$. Note that $\Omega$ is bounded if and only if $\widetilde\Omega$
is bounded, by the compactness of the fibers. Let $f$ be the first eigenfunction of the problem
$\triangle^Nu+\lambda u=0$ in $\Omega$ with Dirichlet boundary conditions, that can be assumed to be
positive in $\Omega$.

Set $X=-\grad^N\big(\log f\big)$, so that $\mathrm{div}^N(X)-\vert X\vert^2=\lambda_1(\Omega)$
is constant in $\Omega$.
If $\widetilde X$ is the horizontal lifting of $X$, then clearly $\vert\widetilde X_q\vert=
\vert X_{\pi(q)}\vert$ for all $q\in\widetilde\Omega$. Moreover, by Lemma~\ref{thm:divergence},
since $H=0$, $\mathrm{div}^M(\widetilde X)_q=\mathrm{div}^N(X)_{\pi(q)}$. Using \eqref{eq:estimatefundtone},
we then obtain:
\[\lambda^*(\widetilde\Omega)\ge\inf_{\widetilde\Omega}\big[\mathrm{div}^M(\widetilde X)-\vert\widetilde X\vert^2\big]=
 \inf_{\Omega}\big[\mathrm{div}^N(X)-\vert X\vert^2\big]=\lambda^*(\Omega).\] This proves Lemma \eqref{thm:equalfundtone}.
\end{proof}
\begin{cor}
Assume that the fibers of $\pi$ are compact and minimal. Then, $\sigma_{\mathrm{ess}}(M)=\emptyset$
if and only if $\sigma_{\mathrm{ess}}(N)=\emptyset$.
\end{cor}
The above result applies in particular to
Riemannian coverings, yielding the following
\begin{cor}
If $M$ is a finite covering of $N$, then $\sigma_{\mathrm{ess}}(M)\ne\emptyset$ if and only if
$\sigma_{\mathrm{ess}}(N)\ne\emptyset$.
\end{cor}
\subsection{Proof of Theorem  \ref{mainThm}} The item ii. of Theorem \ref{mainThm} follows from this Lemma \ref{thm:equalfundtone}.
For if we take a sequence of compact sets $K_{1}\subset K_{2} \subset\cdots$ with $N=\cup_{i=1}^{\infty}K_{i}$. Likewise we have $M=\cup_{i=1}^{\infty}\widetilde{K}_{i}$, where $\widetilde{K}_{i}=\pi^{-1}(K_{i})$. By the proof of \eqref{thm:charinfessspectrum} we have that  $\inf\sigma_{\mathrm{ess}}(N)=\lim_{i\to \infty}\lambda^{\ast}(N\setminus K_{i})$ and $\inf\sigma_{\mathrm{ess}}(M)=\lim_{i\to \infty}\lambda^{\ast}(M\setminus \widetilde{K}_{i})$. However,  $\lambda^{\ast}(N\setminus K_{i})=\lambda^{\ast}(M\setminus \widetilde{K}_{i})$,   by the Lemma \eqref{thm:equalfundtone}.
Before we prove item i. we need the following lemma.
\begin{lem} Let $\pi\colon M\to N$ be a Riemannian submersion with compact minimal fibers $\mathcal{F}$.  If $f\in L^{2}(N)$ and $\triangle^{N} f \in L^{2}(N)$ then $\widetilde{f}\in L^{2}(M)$ and $\triangle^{M} \widetilde{f} =\widetilde{\triangle^{N} f}\in L^{2}(M)$. In other words, if $f\in \mathrm{Dom}(\triangle^{N})$ then $\widetilde{f}\in \mathrm{Dom}(\triangle^{M})$.\label{fabio}
\end{lem}
\begin{proof} Let $\widetilde{f}=f \circ \pi$ be the lifting of $f$. By Fubini's Theorem we have
$$\int_{M}\widetilde{f}^{2}dM=\int_{N}\big(\int_{ \mathcal{F}_{p}}f^{2} \,d\mathcal{F}_{p}\big)dN=\mathrm{vol}(\mathcal{F}_{p})\int_{N}f^{2}dN<\infty.$$
This shows that $\widetilde{f}\in L^{2}(M)$. To  show that $\triangle^{M} \widetilde{f}\in L^{2}(M)$ we will show that $\triangle^{M} \widetilde{f} =\widetilde{\triangle^{N} f}$.  Every $\varphi \in C^{\infty}_{0}(M)$ can be decomposed as $\varphi =\varphi_{1}+\varphi_{2}$ where $\varphi_{1}$ is constant along the fibers $\mathcal{F}$  and $(\varphi_{2})_{\mathrm{av}}=0$, see \cite{bordoni}. Moreover, $\varphi_{1}$ and $\varphi_{2}$ has compact support. Observe that we can define $\psi \colon N\to\mathbb{R}$ by $\psi (\pi (p))=\varphi_{1} (p)$ so that $\varphi_{1}=\widetilde{\psi}$. By the Lemma \ref{thm:laplacianlift} we have that
$\triangle^M \varphi_{1}(p)=\triangle^N\psi(\pi(p))$ for every $p\in M$. By Corollary \ref{thmmediazero} $\int_{M}\widetilde{f}\,\triangle^{M} \varphi_{2} dM=0$, therefore
\begin{eqnarray} \int_{M}\widetilde{f}\,\triangle^{M} \varphi \,dM & = & \int_{M}\widetilde{f}\,\triangle^{M} \varphi_{1} dM \nonumber \\
& =&  \!\int_{N}\!\big(\int_{\mathcal{F}_{p}}\!f\,\triangle^{M} \varphi_{1}d\mathcal{F}_{p}\big)dN \nonumber \\
&=& \!\int_{N}\!\big(\!f\,\triangle^{N} \psi  \int_{\mathcal{F}_{p}}\,d\mathcal{F}_{p}\big)dN \nonumber\\
 &=& \! \mathrm{vol}(\mathcal{F}_{p})\int_{N}\!\!f\,\triangle^{N} \psi  dN \nonumber \\
 &=& \! \mathrm{vol}(\mathcal{F}_{p})\int_{N}\!\!\psi \,\triangle^{N} f  dN \nonumber \\
 &=& \!\int_{N}\!\big(\! \int_{\mathcal{F}_{p}}\psi \,\triangle^{N} f \,d\mathcal{F}_{p}\big)dN \nonumber \\
 &=& \int_{M}\widetilde{\psi \,\triangle^{N} f} dM \nonumber \\
 &=& \int_{M}\varphi_{1}\widetilde{\triangle^{N} f} dM \nonumber\\
 &=&  \int_{M}\varphi\widetilde{\triangle^{N} f} dM \nonumber \end{eqnarray}\end{proof}
 To  show that $\sigma_{p}(N)\subset \sigma_{p}(M)$ we take $\lambda \in \sigma_{p}(N)$ and  $f\! \in \! L^{2}(N)$ with $-\triangle^{N}f=\lambda f$ in distributional sense. This implies that $\widetilde{-\triangle^{N}f}=\lambda \widetilde{f}$.  By Lemma \ref{fabio}, $-\triangle^{M}\widetilde{f}=\lambda \widetilde{f}$ showing that $\lambda\in \sigma_{p}(M)$.
To show that $\sigma_{\mathrm{ess}}(N)\subset \sigma_{\mathrm{ess}}(M)$ we take $\mu \in \sigma_{\mathrm{ess}}(N)$. There exists an orthonormal sequence of  functions $f_{k}\in \mathrm{Dom} (\triangle)$ such that $\Vert-\triangle^{N} f_{k}-\mu f_{k}\Vert_{L^{2}(N)}\to 0$ as $k\to \infty$. By Lemma \ref{fabio}, we have that $\widetilde{f_{k}}\in \mathrm{Dom}(\triangle^{M})$. Now
\begin{eqnarray}\Vert-\triangle^{M} \widetilde{f_{k}}-\mu \widetilde{f_{k}}\Vert_{L^{2}(M)}^{2}&=& \int_{M}\vert -\triangle^{M} \widetilde{f_{k}}-\mu \widetilde{f_{k}}\vert^{2}dM\nonumber \\
&=& \int_{N}\int_{\mathcal{F}_{q}}\vert-\triangle^{N}f_{k}-\mu f_{k}\vert^{2}d\mathcal{F}_{q}dN \nonumber \\
&=& \mathrm{vol}(\mathcal{F}_{q})\int_{N}\vert -\triangle^{N}f_{k}-\mu f_{k}\vert^{2}dN\nonumber \\
&=&\mathrm{vol}(\mathcal{F}_{q})\,\Vert-\triangle^{N} f_{k}-\mu f_{k}\,\Vert_{L^{2}(N)}^{2}\to 0\nonumber
\end{eqnarray}This shows that $\mu \in \sigma_{\mathrm{ess}}(M)$,  the proof of Theorem~\ref{mainThm} is concluded.
\begin{cor}\label{thm:corGK}
Let $G$ be a Lie group endowed with a bi-invariant metric. Then, $\sigma_{\mathrm{ess}}(G)$ is empty if and only
if for some (hence for any) compact subgroup $K\subset G$, the Riemannian homogeneous space $G/K$ has
empty essential spectrum.
\end{cor}
\begin{proof}
Apply Theorem~\ref{mainThm} to the Riemannian submersion $G\mapsto G/K$, which has minimal and compact fibers.
\end{proof}
Other interesting examples of applications of Theorem~\ref{mainThm} arise from non compact Lie groups.
\begin{example-n}
Consider the $2\times 2$ special linear group $\mathrm{SL}(2,R)$. There exists a $2$-parameter family of
left-invariant Riemannian metrics $g_{\kappa,\tau}$, with $\kappa<0$ and $\tau\ne0$, for which
$\big(\mathrm{SL}(2,\R),g_{\kappa,\tau}\big)\to \mathbb{H}^{2}(\kappa)$ is a Riemannian submersion
with geodesic fibers diffeomorphic to the circle $\mathds S^1$. An explicit description
of these metrics can be found, for instance, in \cite{torralbo}.
Endowed with these metrics, $\mathrm{SL}(2,\R)$ is one of the eight homogeneous Riemannian
$3$-geometries, as classified in \cite{Sco}, and its isometry group has dimension $4$.
\begin{prop}\label{thm:SL2R}
For all $\kappa<0$ and $\tau\ne0$, \[\sigma\big(\mathrm{SL}(2,\R),g_{\kappa,\tau}\big)=\sigma_{\mathrm{ess}}
\big(\mathrm{SL}(2,\R),g_{\kappa,\tau}\big)=\left[\displaystyle-\tfrac\kappa4,+\infty\right).\]
\end{prop}
\begin{proof}
It is   known that the spectrum $\sigma\big(\mathbb{ H}(\kappa)\big) =\sigma_{\mathrm{ess}}
\big(\mathds H(\kappa)\big)=\left[-\tfrac\kappa4,+\infty\right)$, see \cite{donnelly}.
By Lemma~\ref{thm:equalfundtone}
\[\lambda^*\big(\mathrm{SL}(2,\R),g_{\kappa,\tau}\big)=\lambda^*\big(\mathds H(\kappa)\big)=-\tfrac\kappa4,\]
hence $\sigma\big(\mathrm{SL}(2,\R),g_{\kappa,\tau}\big)\subset\left[-\tfrac\kappa4,+\infty\right)$.
On the other hand, by Theorem~\ref{mainThm}
\[\left[-\tfrac\kappa4,+\infty\right)=\sigma_{\mathrm{ess}}\big(\mathds H(\kappa)\big)\subset
\sigma_{\mathrm{ess}}\big(\mathrm{SL}(2,\R),g_{\kappa,\tau}\big).\] This proves the proposition.
\end{proof}
\end{example-n}
\section{Mean curvature of geodesic spheres versus  mean curvature of the fibers. Proof of Theorem~\ref{mainThm2}.}
\label{sec:mencurvfibersvscurvspheres}
We will now drop the minimality and the compactness assumption on the fibers, however, we will make some
 assumptions on the curvature of the base and the fibers of the submersion.
Assume that $(N,g^N)$ has a \emph{pole} $p_0$ or more generally has a point $p_{0}$ with {\em thin cut locus}, see the Introduction.
For $p\in N\setminus\{p_0\}$, let $\gamma_p\colon[0,1]\to \!N$ be the unique affinely parameterized
geodesic in $(N,g^N)$ such that $\gamma_p(0)=p_0$ and $\gamma_p(1)=p$.
The \emph{radial curvature} function of $(N,g^N\!\!, p_0)$, denoted
by $\kappa_{p_0}:N\to\R$, is defined by $\kappa_{p_0}(p)=\max\limits_\sigma\mathrm{sec}(\sigma)$, where
$\mathrm{sec}$ is the section curvature and the maximum is taken over all $2$-planes $\sigma\subset T_pN$
containing the direction $\gamma'_p(1)$. Finally, let us denote by $\rho_{p_0}:N\to[0,+\infty)$
the distance function in $N$ given by $\rho_{p_0}(p)=\mathrm{dist}_N(p,p_0)$.

We are now ready for

\begin{proof}[Proof of Theorem~\ref{mainThm2}] Assume first  that
 $\pi:M\to N$ is a Riemannian submersion satisfying the following assumptions:
\begin{itemize}
\item[(a)] $(N,g^N)$ has a pole $p_0$.
\item[(b)] the function
$h(q)=(\triangle^N\rho_{p_0})_{\pi(q)}+g^N\big((\grad^N\rho_{p_0})_{\pi(q)},\mathrm d\pi_q(H_q)$
is \emph{proper}.
\end{itemize}
Consider the lifting  $\widetilde\rho_{p_0}\colon \!M\!\to\!\R$ defined by $\widetilde\rho_{p_0}=\rho_{p_0}\circ\pi$. Then
by Lemma~\ref{thm:laplacianlift} we have that $h=\triangle^M\widetilde\rho_{p_0}$. Moreover, by \eqref{eq:gradlifting},
\[\vert\mathrm{grad}^M\widetilde\rho_{p_0}\vert=\vert\mathrm{grad}^N\rho_{p_0}\vert\equiv1.\]
If $K_1\subset K_2\subset\cdots\subset K_n\subset\ldots$ is a compact exhaustion of $M$,
then by \eqref{eq:secondestimate} applied to $X=\grad^M\widetilde\rho_{p_0}$
\begin{equation}\lambda^*(M\setminus K_n)\ge\frac{1}{4}\big[\inf_{M\setminus K_n}h\big]^2.\label{eqFinal}\end{equation}
Since $h$ is proper, the right-hand side in the above inequality tends to $+\infty$
as $n\to\infty$, thus, by Proposition~\ref{thm:charinfessspectrum}, $\sigma_{\mathrm{ess}}(M)=\emptyset$.

 If $p_{0}$ has thin cut locus, the same proof above holds, since
 $X=\grad^M\widetilde\rho_{p_0}$ satisfies the Proposition \ref{thm:est1} and therefore \ref{eqFinal},
 see Remark \ref{thincutlocus}.
\end{proof}

\begin{cor}\label{thm:cormainthm2}
If the fibers of $\pi\colon M \to N$ are compact and if  the function $l(p)=\triangle^N\rho_{p_0}(p)-\max\limits_{q\in\mathcal F_p}\Vert H_q\Vert$
is proper then $\sigma_{ess}(M)=\emptyset$.
\end{cor}
A different statement can be obtained in terms of radial curvature.
\begin{cor}\label{thm:radialcurvature}
Assume that $G:\left[0,+\infty\right[\to\R$ is a smooth function such that:
\begin{equation}\label{eq:disradialcurvature}
\kappa_{p_0}(p)\le-G\big(\rho_{p_0}(p)\big)
\end{equation}
for all $p\in N$. Denote by $\psi:\left[0,+\infty\right[\to\R$ the solution of the Cauchy problem:
\[\psi''(t)=G(t)\psi(t),\quad\psi(0)=0,\quad\psi'(0)=1,\]
and set $\ell(t)=(n-1)\psi'(t)/\psi(t)$, $t>0$.
If
\begin{equation}\lim\limits_{q\to\infty}\Big[\ell\big(\rho_{p_0}\big(\pi(q)\big)+
g^N\big((\nabla^N\rho_{p_0})_{\pi(q)},\mathrm d\pi_q(H_q)\big)\Big]=+\infty\label{b''}\end{equation}
then $\sigma_{ess}(M)=\emptyset$.
\end{cor}
\begin{proof}
Using the Hessian Comparison Theorem \cite[Chapter~2]{GreeneWu}, under the assumption
\eqref{eq:disradialcurvature} one has:
\begin{equation}\label{eq:hessiancomparison}
\mathrm{Hess}(\rho^N)\ge\frac{\psi'}\psi\cdot\big(g^N-\mathrm d\rho_{p_0}\otimes\mathrm d\rho_{p_0}\big).
\end{equation}
Considering an orthogonal basis of $T_pN$ of the form $\{\nabla^N\rho_{p_0},e_1,\ldots,e_{n-1}\}$, where
$\{e_1,\ldots,e_{n-1}\}$ is an orthonormal basis of $T_p\mathcal S_p$, and taking the trace of the symmetric bilinear
forms in the two sides of \eqref{eq:hessiancomparison}, we get
\begin{equation}\label{eq:laplaciancomparison}
\triangle^N\rho_{p_0}\ge(n-1)\frac{\psi'}\psi.
\end{equation}
It is clear that (\ref{b''}) implies that   $h(q)=\triangle^N\rho_{p_0}+g^N\big(\grad^N\rho_{p_0},\mathrm d\pi_q(H_q))$ is proper.
\end{proof}
\subsection{Warped products}
Let $(N,g^N)$ and $(F,g^F)$ be Riemannian manifolds and let $\psi:N\to\R^+$ be a smooth
function. The \emph{warped product manifold} $M=N\times_\psi F$ is the product manifold $N\times F$
endowed with the Riemannian metric $g^N+\psi^2g^F$. It is immediate to see that the projection
$\pi:M\to N$ onto the first factor is a Riemannian submersion, with fiber $\mathcal F_p=\{p\}\times F$.
Among Riemannian submersions, warped products are characterized by the following properties:
\begin{itemize}
\item the horizontal distribution is integrable, and its leaves are totally geodesic;
\item the fibers are totally umbilical.
\end{itemize}
For warped products, the results of the paper can be stated in a more explicit form in terms
of the warping function $f$. The mean curvature of the fibers are given by
\begin{equation}\label{eq:meancurvaturefiberwarped}
H=-\mathrm{dim}(F)\,\frac{\grad^M\widetilde\psi}{\widetilde\psi},
\end{equation}
where $\widetilde{\psi}$ is the lifting of $\psi$.

\begin{prop}
Let $M=N\times_\psi F$ be a warped product, with $F$ compact.
\begin{itemize}
\item[(a)] If $\sigma_{\mathrm{ess}}(N)\ne\emptyset$, and
$0<\inf_N\psi\leq \sup_N\psi<+\infty$.

\noindent Then $\sigma_{\mathrm{ess}}(M)\ne\emptyset$.
\smallskip

\item[(b)] If $(N,g^N)$ has a pole $p_0$ and the function
\[\triangle^N\rho_{p_0}-\frac1\psi g^N\big(\nabla^N\rho_{p_0},\nabla^N\psi\big)\] is proper, then
$\sigma_{\mathrm{ess}}(M)\ne\emptyset$.
\end{itemize}
\end{prop}
\begin{proof}
Part (a) follows from Proposition~\ref{lemma3.5}, observing that the volume
of the fiber $\mathcal F_p=\{p\}\times F$ equals $\psi(p)^{\mathrm{dim}(F)}\mathrm{vol}(F)$.

Part (b) follows from Theorem \ref{mainThm2} and formula \eqref{eq:meancurvaturefiberwarped}.
\end{proof}

\subsection{Example}\label{examples}
Let $\mathbb{R}^{n}=[0, \infty)\times \mathbb{S}^{n-1}/\sim$ endowed with the smooth  metric  $ds^{2}=dt^{2}+ f^{2}(t)d\theta^{2}$, where $f(0)=0$, $f'(0)=1$. The equivalence relation $\sim$ is the following $$(t, \theta)\sim (s, \alpha) \Leftrightarrow t=s=0\,\, \mathrm{or}\,\,t=s>0\,\, \mathrm{and}\,\, \theta=\alpha.$$ The radial sectional curvatures $K^{\mathrm{ rad}}$ along the geodesic issuing the origin $0=\{0\}\times \mathbb{S}^{n-1}/\sim$ is given by $K^{\mathrm{ rad}}(t, \theta)=-\displaystyle\frac{f''(t)}{f(t)}$. Let us consider $W=\mathbb{R}^{n}\times \mathbb{S}^{1}$  with metric $dw^{2}=ds^{2}+ g^{2}(\rho(x))d^{2}_{\mathbb{S}^{1}}$, where $\rho$  is the distance function to the origin in $(\mathbb{R}^{n}, ds^{2})$ and $g\colon [0, \infty)\to (0, \infty) $ is a smooth function.
Choosing $f(t)= te^{t^{2}}~$ we have   $K^{\mathrm{ rad}}(t, \theta)=-4t^{3}-6t$, thus
$\lim_{t\to\infty}K^{\mathrm{ rad}}(t, \theta)=-\infty$. By Donnelly-Li's Theorem  \cite{donnelly-li}, the spectrum of $(\mathbb{R}^{n}, ds^{2})$ is discrete. Choosing
$g(t)=e^{t-t^2}$. An easy computation yields that
\begin{itemize}\item[i.] The volume $\mathrm{vol}(W, dw^2)=\infty.$
\item[ii.] The limit $\mu= \limsup_{r\to \infty} \displaystyle\frac{\log (\mathrm{vol}(B_W(r)))}{r}< \infty$, where $B_W(r)$ is the geodesic ball centered at a point $p=(0, \xi)\in W$ and radius $r$.
\end{itemize}The items i.\!\! and ii. imply by Brooks' Theorem \cite{brooks}  that $\sigma_{\mathrm{ess}}(W)\neq \emptyset$. This gives an example of a Riemannian submersion $\pi \colon (W,dw^2) \to (\mathbb{R}^{n}, ds^2)$ where the base space is discrete but the total space is not, while the fiber is compact but not minimal.

 An example of a Riemannian submersion $\pi \colon (\mathbb{R}^{n}\times \mathbb{S}^{1},dw^2) \to (\mathbb{R}^{n}, ds^2)$ where the total space is discrete but the base space is not, while the fiber is compact but not minimal is presented in \cite[Proposition 4.3]{baider}.
\end{section}


\begin{thebibliography}{10}

\bibitem{baider} \textsc{A. Baider}, \emph{Noncompact Riemannian manifolds with discrete spectra}, J. Diff.\ Geom.\
\textbf{14} (1979), 41--57.

\bibitem{bessa-montenegro} \textsc{G. P. Bessa, J. F. Montenegro}, \emph{Eigenvalue estimates for submanifolds
with locally  bounded mean  curvature}, Ann.\ Global Anal.\ Geom.\ \textbf{24} (2003), 279--290.

\bibitem{bessa-montenegro2} \textsc{G. P. Bessa, J. F. Montenegro}, \emph{An extension of Barta's
theorem and geometric applications}, Ann.\ Global Anal.\ Geom.\  \textbf{31}  (2007),  no.\ 4, 345--362.

\bibitem{bessa-jorge-montenegro} \textsc{G. P. Bessa, L. Jorge, J. F. Montenegro}, \emph{The  spectrum of the Martin-Morales-Nadirashvili minimal surfaces is discrete.} to appear in J. Geom.\ Anal.\

\bibitem{bordoni}
\textsc{M. Bordoni}, \emph{Spectral estimates for submersions with fibers of basic mean curvature},
An.\ Univ.\ Vest Timi\c s.\ Ser.\ Mat.-Inform.\  \textbf{44}  (2006),  no.\ 1, 23--36.

\bibitem{brooks}
\textsc{R. Brooks}, \emph{A relation between growth and the spectrum of the Laplacian},
Math.\ Z.\  \textbf{178}  (1981),   501--508.

\bibitem{davies} \textsc{E. B. Davies}, \emph{Spectral theory and differential operators}, Cambrigde University Press, 1995.

\bibitem{donnelly} \textsc{H. Donnelly}, \emph{On the essential spectrum of a complete Riemannian manifold},
Topology \textbf{20}  (1981), no.\ 1, 1--14.

\bibitem{donnelly1}
\textsc{H. Donnelly}, \emph{Negative curvature and embedded  eigenvalues.} Math. Z. \textbf{203}, (1990), 301--308.

\bibitem{donnelly-garofalo}
\textsc{H. Donnelly, N. Garofalo}, \emph{Riemannian manifolds whose Laplacian have purely continuous spectrum.} Math. Ann. \textbf{293}, (1992), 143--161.

\bibitem{donnelly-li} \textsc{H. Donnelly,  P. Li}, \emph{Pure point spectrum and negative curvature for
noncompact manifolds}, Duke Math.\ J. \textbf{46} (1979), 497--503.

\bibitem{escobar}
\textsc{E. Escobar}, \emph{On the spectrum of the Laplacian on complete Riemannian manifolds.} Comm. Partial Differ. Equations \textbf{11}, (1985), 63--85.

\bibitem{Gra67}
\textsc{A. Gray}, \emph{Pseudo-{R}iemannian almost product manifolds and
  submersions}, J. Math.\ Mech.\ \textbf{16} (1967), 715--737.

\bibitem{GreeneWu} \textsc{R. E. Greene, H. Wu}, \emph{Function theory on manifolds which possess a pole},
Lecture Notes in Mathematics, 699. Springer, Berlin, 1979.

\bibitem{harmer} \textsc{M. Harmer}, \emph{Discreteness of the spectrum of the Laplaceian and Stochastic incompleteness}, J. Geom.\ Anal.\ \textbf{19} (2009), 358--372.

\bibitem{karp}
\textsc{L. Karp}, \emph{Noncompact manifolds with purely continuous spectrum}, Mich. Math. J. \textbf{31}, (1984), 339--347.

\bibitem{kleine1} \textsc{R. Kleine}, \emph{Discreteness conditions for the Laplacian on complete noncompact Riemannian manifolds}, Math.\ Z. \textbf{198} (1988), 127--141.

\bibitem{kleine2} \textsc{R. Kleine}, \emph{Warped products with discrete spectra}, Results Math.\ \textbf{15} (1989), 81--103.

\bibitem{One66}
\textsc{B. O'Neill}, \emph{The fundamental equations of a submersion}, Michigan
  Math.\ J.\ \textbf{13} (1966), 459--469.

\bibitem{One67}
\leavevmode\vrule height 2pt depth -1.6pt width 23pt, \emph{Submersions and
  geodesics}, Duke Math.\ J.\ \textbf{34} (1967), 363--373.

%\bibitem{One83}
%\leavevmode\vrule height 2pt depth -1.6pt width 23pt, \emph{Semi-{R}iemannian
%  geometry}, vol.\ 103 of Pure and Applied Mathematics, Academic Press Inc.
%  [Harcourt Brace Jovanovich Publishers], New York, 1983.
%\newblock With applications to relativity.

\bibitem{rellich} textsc{F. Rellich}, \emph{\"{U}ber das asymptotische Verhalten der L\"{o}sungen von $\triangle u+\lambda u=0$ in unendlichen Gebieten.} Jahresber. Dtsch. Math.-Ver. \textbf{53}, (1943), 57--65.

\bibitem{Sco} \textsc{P. Scott}, \emph{The geometries of 3-manifolds}, Bull.\ London Math.\ Soc.\ \textbf{15}
(5) (1983), 401--487.

\bibitem{tayoshi} \textsc{T. Tayoshi}, \emph{On the spectrum of the Laplace-Beltrami operator on noncompact surface},
Proc.\ Japan Acad.\ \textbf{47}, (1971), 579--585.

\bibitem{torralbo} \textsc{F. Torralbo}, \emph{Rotationally invariant constant mean curvature surfaces in
homogeneus $3$-manifolds}, preprint 2009, \texttt{arXiv:0911.5128v1}.

\end{thebibliography}
\end{document}